\documentclass[submission,copyright,creativecommons]{eptcs}
\providecommand{\event}{GASCom 2026} 

\usepackage{iftex}

\ifpdf
  \usepackage{underscore}         
  \usepackage[T1]{fontenc}        
\else
  \usepackage{breakurl}           
\fi

\newtheorem{theorem}{Theorem}
\newtheorem{lemma}[theorem]{Lemma}
\newtheorem{conjecture}[theorem]{Conjecture}
\newtheorem{example}[theorem]{Example}
\newtheorem{proposition}[theorem]{Proposition}
\newtheorem{definition}[theorem]{Definition}
\newtheorem{remark}[theorem]{Remark}

\usepackage{lipsum}
\usepackage{xcolor}
\usepackage{algorithm}
\usepackage{algpseudocode}
\usepackage{float} 

\newcommand{\s}{\sigma}
\newcommand{\q}{\tau}

\newcommand{\1}{^{-1}}
\renewcommand{\b}{\beta}
\renewcommand{\l}{\ell}
\renewcommand{\c}{\cdots}
\newcommand{\jo}{\vee_R}
\newcommand{\w}{\leq_R}
\newcommand{\Inv}{\mathrm{Inv}}
\renewcommand{\t}[2]{\left(
        #1\ #2
    \right)}

\renewcommand{\r}{\rho}
\renewcommand{\o}[1]{\noverline{#1}}
\newcommand{\Neg}{\mathrm{Neg}}
\newcommand{\Nsp}{\mathrm{Nsp}}
\newcommand{\bncut}[1]{(B_{#1},B^*_{#1})}

\usepackage{mathtools}
\newcommand\nindent{.5pt}
\newcommand\noverline[1]{%
  \kern\nindent\overline{\kern-\nindent#1\kern-\nindent}\kern\nindent}

\title{Computing Joins in the Weak Order of Type $B$ Coxeter Groups: an Algorithmic Approach}
\author{Riccardo Biagioli
\institute{Dipartimento di Matematica\\ Università di Bologna, Italy}
\email{riccardo.biagioli2@unibo.it}
\and
Lorenzo Perrone
\institute{Dipartimento di Matematica\\ Università di Bologna, Italy}
\email{lorenzo.perrone8@unibo.it}
}

\def\titlerunning{Joins in Type $B$ Weak Order}
\def\authorrunning{R. Biagioli \& L. Perrone}

\begin{document}
\maketitle

\begin{abstract}
We present an algorithm for computing the join of two elements in the weak order of the Coxeter group of type $B$. This extends Markowsky's algorithm for computing joins of standard permutations to signed permutations, and allows us to confirm a conjecture of Dyer concerning a geometric interpretation of these joins.
\end{abstract}

\section{Introduction}

In a very influential paper, Dyer~\cite{dyer2019} proposed several conjectures regarding the weak order of Coxeter groups. One of these conjectures, which will be discussed in detail in the following sections, has been  the primary motivation for the present study.

\begin{figure}[h]
    \centering
    \includegraphics[scale=1]{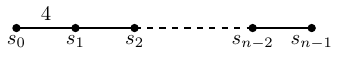}
    \caption{The Coxeter graph of type $B_n$.}
    \label{fig:Bngraph}
\end{figure}

It is well known that the Coxeter system of type $B_n$, whose Coxeter graph is displayed in Figure \ref{fig:Bngraph} admits a combinatorial interpretation as the group $S^B_n$ of signed permutations, namely the set of bijections $\s:[\pm n]\rightarrow [\pm n]$ satisfying the symmetric property $\s(-a)=-\s(a)$ for every $a\in[n]$, where we denote $[n]=\{1,\ldots,n\}$ and $[\pm n]=\{-n,\ldots,-1,1,\ldots n\}.$
We refer to \cite[\S8.1]{bjorner2005combinatorics} for any undefined notation.

To avoid the use of minus signs, we write $\o{i}$ instead of $-i$. Let $i\neq\pm j\in[\pm n]$; we denote by $\t{i}{j}\in S^B_n$ the \emph{transposition} that swaps $i$ with $j$ and $\o{i}$ with $\o{j}$, and by $\t{i}{\o{i}}$ the one that exchanges $i$ with $\o{i}$. Note that, naturally, $\t{i}{j}$, $\t{j}{i}$, $\t{\o{i}}{\o{j}}$ and $\t{\o{j}}{\o{i}}$ all represent the same signed permutation. 
Through this combinatorial interpretation, we take as set of Coxeter generators for $S^B_n$ the set $S_B=\{s_0\!=\!\t{1}{\o{1}},s_1\!=\!\t{1}{2},\ldots,s_{n-1}\!=\!\t{n-1}{n}\}$.

In what follows, we will use both the \emph{complete notation} of a signed permutation $$\s=\s(\o{n})\c\s(\o{1})\s(1)\c\s(n)$$ and its \emph{window notation} $\s=\s(1)\s(2)\cdots\s(n)$ which exploits the symmetry $\s(\o{i})=\o{\s(i)}$. We will specify which notation we use only when it is strictly necessary.

Recall that, given a Coxeter system $(W,S)$, its set of reflections is $T=\bigcup_{w\in W}wSw\1$ and the {\em left-reflection set} of an element $w\in W$ is $T_L(w)=\{t\in T\mid \l(tw)<\l(w)\}$.
The \emph{(right) weak order} on $(W,S)$ is a partial order defined by the prefix property: for any $u,v\in W$, we have $u\w v$ if and only if a reduced expression for $u$ is a prefix of a reduced expression for $v$. It can also be characterized in terms of left-reflection sets, namely 
\begin{equation}
    u\w v\ \text{ if and only if }\ T_L(u)\subseteq T_L(v).
\end{equation}
The poset $(W,\w)$ is a meet-semilattice, and in particular it is a lattice whenever the group $W$ is finite. On the other hand, it is never a lattice when $W$ is infinite.
For further details on $(W,\w)$, we refer the reader to \cite[Chapter~3]{bjorner2005combinatorics}.

In the combinatorial description $(S^B_n,S_B)$, the set of reflections corresponds to the set of transpositions $T=\{\t{i}{j}\mid i\neq j \in [\pm n]\}$ introduced above. Moreover, in this interpretation, it is useful to introduce the following well known statistics of a signed permutation $\s\in S^B_n$:
\begin{itemize}
    \item $\Inv(\s)=\{(i,j)\in [n]\times[n]\mid i<j, \s(i)>\s(j)\}$, the set of \emph{inversions} of $\s$;
    \item $\Neg(\s)=\{i\in [n]\mid \s(i)<0\}$, the set of \emph{negative entries} of $\s$;
    \item $\Nsp(\s)=\left\{\{i,j\}\subseteq [n]\mid i\neq j, \s(i)+\s(j)<0 \right\}$, the set of \emph{negative sum pairs} of $\s$.
\end{itemize}
The computation of these statistics provides a direct formula for the Coxeter length:
$$\l(\s)=|\Inv(\s)|+|\Neg(\s)|+|\Nsp(\s)|.$$
This identity reflects a decomposition of the left-reflection set of $\sigma$ into three disjoint subsets:
\begin{equation}\label{eq:leftref}
    T_L(\s) = \left\{ \t{i}{j} \mid  (i,j) \in \Inv(\s\1)  \right\} \sqcup \left\{ \t{i}{\o{i}} \mid  i \in  \Neg(\s\1) \right\} \sqcup \left\{ \t{i}{\o{j}} \mid  \{i,j\} \in \Nsp(\s\1) \right\}.
\end{equation}

A proof of last equality can be found in \cite[Chapter~8]{bjorner2005combinatorics} and \cite{Yu2024}.
It is important to observe that the left reflection set $T_L(w)$ uniquely identifies the element $w$. Therefore, to determine the join of two elements in $(W,\leq_R)$, it suffices to compute its left reflection set. For instance, Theorem 10-3.25 of \cite{Reading2016} gives two descriptions of the left reflection set of the join of two elements in terms of closure operations. As we will see below, the conjecture that motivated our work is formulated in a similar way.

In \cite{dyer2019}, Dyer introduced the extended weak order, a partial order generalizing the weak order, and conjectured that it is a lattice even when the underlying Coxeter group is infinite. In the same paper, he formulated a conjectural description of the join of two elements in this poset. It is not yet known whether this conjecture holds even for finite Coxeter systems. It has been proved for Coxeter groups of types $I$ and $A$ in~\cite{BiagioliPerroneA}, and verified for the exceptional types $F_4$, $H_3$, and $H_4$ using SageMath~\cite{sagemath}. Establishing this conjecture for Coxeter systems of type $B$ was the main motivation for the present paper. For the original statement, we refer the reader to [4, §2.8]. Here, we give an equivalent formulation for finite Coxeter groups in terms of the {\em Bruhat preclosure} recently introduced by Dermenjian in~\cite[\S3]{dermenjian2025}.

\smallskip

For each Coxeter system $(W,S)$, we recall that the Bruhat graph $B(W)$ is the directed graph having $W$ as vertex set, and where for any $u,v\in W$ there is an arrow $u\xrightarrow{t}v$ whenever $v=tu$ and $\l(u)<\l(v)$.

\begin{definition}[$A$-Bruhat path]
    Let $W$ be a Coxeter group; consider the Bruhat graph $B(W)$ associated with $W$ and $A\subseteq T$. An $A$-\emph{Bruhat path} 
    is any directed path in $B(W)$ starting from the identity element $e$ of $W$, whose edges are labeled by elements in the set $A$. 
    If a set $A\subseteq T$ can be written as a union of two left-reflection sets, i.e. $A=T_L(u)\cup T_L(v)$ for some $u,v\in W$, then we use the expression $(u,v)$-Bruhat paths instead of $A$-Bruhat paths.
\end{definition}

\begin{definition}[Bruhat preclosure]
    Given a set of reflection $A\subseteq T$; its Bruhat preclosure is 
    $$\overline{A}=\{t\in T\mid \text{there is an $A$-Bruhat path from }e \text{ to }t\}.$$
If $A=\overline{A}$ we say it is \emph{preclosed}.
\end{definition}

Observe that this operation is not a closure operator: in general, 
$\overline{\overline{A}}\neq \overline{A}$. Indeed, counterexamples were found by Dermenjian in Coxeter groups of types $H_3$ and $F_4$; see~\cite[\S3]{dermenjian2025}. 
Consequently, the following conjecture is of a rather different nature from Theorem~10-3.25, where the operations involved are proper closure operators.

\begin{conjecture}\label{conjecture}
    Let $W$ be a finite Coxeter group and $u,v\in W$. Then
        $$T_L(u\jo v)=\overline{T_L(u)\cup T_L(v)}.$$
\end{conjecture}

\begin{figure}[h]
    \centering
    \includegraphics[scale=1]{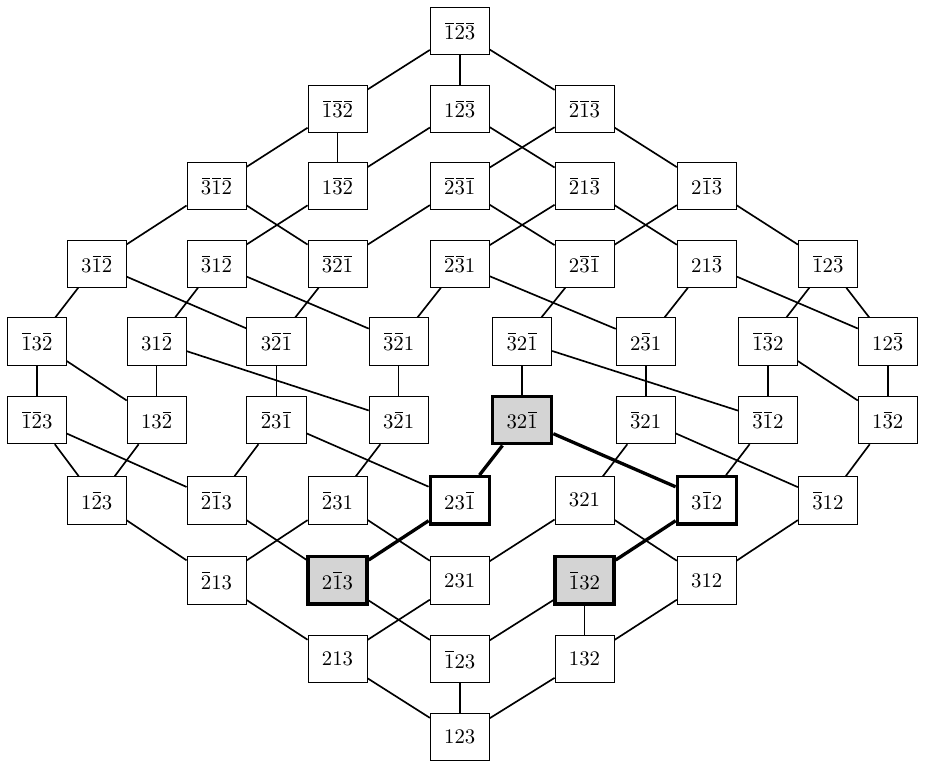}
    \caption{The Hasse diagram of $(S^B_3,\w)$, where signed permutations are in window notation and the join operation of $2\o{1}3$ and $\o{1}32$ is highlighted.}
    \label{fig:Hasse}
\end{figure}

In general, computing the join of two elements in the poset $(W,\w)$ can be quite difficult; see for instance Figure~\ref{fig:Hasse}. 
Conjecture~\ref{conjecture} states that the left-reflection set of the join $u \jo v$ can be determined by considering all reflections reached by $(u,v)$-Bruhat paths. 

In Section~\ref{secalg}, we present the main result in this extended abstract  which is an algorithmic way to compute the join of two signed permutations starting from their complete notation. 
This is a generalization of the algorithm presented in \cite{markowsky1994permutation} for the symmetric group. This algorithm will serve as a key tool in a forthcoming article, where we establish the validity of Dyer’s conjecture for Coxeter groups of type $B$.

Observe that, thanks to equation~\eqref{eq:leftref}, the weak order on $S^B_n$ can be characterized as:
\begin{equation}\label{eq:weakstat}
    \s\w\q \iff \Inv(\s\1)\subseteq\Inv(\q\1),\ \Neg(\s\1)\subseteq\Neg(\q\1),\ \Nsp(\s\1)\subseteq\Nsp(\q\1).
\end{equation}

\section{Algorithm}\label{secalg}

In this section, we describe an algorithm to compute the join of two elements in $(S^B_n,\w)$. Note that when restricted to standard permutations this coincides with the algorithm in \cite[\S 2]{markowsky1994permutation}; therefore it can be considered a generalization of it. 
For simplicity and readability in the definition of the algorithm, we use the complete notation for the two input signed permutations.

\begin{definition}[$n$-cut]
    The {\em $n$-cut} of $\s \in S^B_n$ is the pair $\bncut{\s}$, where $B_\s$ and $B^*_\s$ are the strings made of all entries different from $\o{n}$, respectively, on the left and on the right of the entry $n$ in the complete notation of $\s$.
    More precisely, if $k:=|\s\1(n)|$, we have:
    \begin{enumerate}
        \item if $\s\1(n)<0$, then
         \begin{equation*}
        \begin{split}
            &B_\s=\s(\o{n})\s(\o{n-1})\c\s(\o{k+1}),\\
            &B_\s^*=\s(\o{k-1})\c\s(\o{1})\s(1)\c\s(k-1)\s(k+1)\c\s(n);
            \end{split}
        \end{equation*}
        \item if $\s\1(n)>0$, then
        \begin{equation*}
        \begin{split}              
            &B_\s=\s(\o{n})\s(\o{n-1})\c\s(\o{k+1})\s(\o{k-1})\c\s(\o{1})\s(1)\c\s(k-1),\\
            &B_\s^*=\s(k+1)\c\s(n).
            \end{split}
            \end{equation*}
    \end{enumerate}
    Moreover, in case 1, we write $B^*_\s=B^{*'}_\s B^{*''}_\s$, where
    \begin{equation*}
        B^{*'}_\s=\s(\o{k-1})\c\s(\o{1})\s(1)\c\s(k-1), \qquad B^{*''}_\s=\s(k+1)\c\s(n),
    \end{equation*}
    while in case 2, we write  $B_\s=B'_\s B''_\s$, where 
    \begin{equation*}
        B'_\s=\s(\o{n})\s(\o{n-1})\c\s(\o{k+1}),   \qquad B''_\s=\s(\o{k-1})\c\s(\o{1})\s(1)\c\s(k-1).
    \end{equation*}
    Thus, in case 1, we express the complete notation of the signed permutation as $\s=B_\s nB^{*'}_\s\o{n}B^{*''}_\s$, while in case 2, we can write $\s=B'_\s \o{n}B''_\s nB^*_\s$.
\end{definition}

\begin{example}
    Consider $\s=\o{3}4\o{1}\o{2}21\o{4}3$ and $\q=2\o{1}\o{4}3\o{3}41\o{2}$ in $S^B_4$. Their $4$-cuts are 
    $$\bncut{\s}=(\o{3},\o{1}\o{2}213),\qquad \bncut{\q}=(2\o{1}3\o{3},1\o{2}).$$
    Furthermore, we have $\s\1(4)<0$, so, $B^{*'}_\s=\o{1}\o{2}21$ and $B^{*''}_\s=3$; while $\q\1(4)>0$, so,
    $B'_\q=2\o{1}$ and $B''_\q=3\o{3}$. Hence, $\s=B_\s 4B^{*''}_\s \o{4} B^{*'}_\s$ and $\q=B'_\q \o{4}B''_\q 4 B^*_\q$.
\end{example}

\begin{remark}
    Note that the $n$-cut of a signed permutation $\s\in S^B_n$ determines it uniquely except when $\s\1(n)=\pm 1$, in which case the sign of the entry $\pm n$ is not encoded in the cut. For instance, $\s=4\o{1}3\o{2}$ and $\q=\o{4}\o{1}3\o{2}$ share the same $4$-cut, but they differ in the sign of the entry $4$.
\end{remark}

We define the map $\r:S^B_n\rightarrow S^B_{n-1}$ as $\r(\s):=B_\s B^*_\s\in S^B_{n-1}.$ 
For brevity, we write $\r_\s$ instead of $\r(\s)$. In particular, for $k\in[n-1]$
\begin{equation*}
    \r_\s(k)=\begin{cases}
        \s(k),\quad&\emph{if }0<k<|\s\1(n)|,\\
        \s(k+1),&\emph{if }|\s\1(n)|\leq k\leq n-1.
    \end{cases}
\end{equation*}
To avoid any misunderstanding, by $\r_\s\1$ we always denote the inverse of $\r_\s$ in the group $S^B_{n-1}$; moreover, observe that $\r_\s\1\neq\r_{\s\1}$, for instance, if $\s=4213$, then $\r_\s\1=213\neq 321=\r_{\s\1}.$

The following lemma is fundamental for understanding the map $\r$; specifically, it details how the statistics of a signed permutation are transformed under its application.

\begin{lemma}\label{lem:ncutInvNegNsp}
Let $\s\in S^B_n$; then
    \begin{align*}
        \Inv(\r_\s\1)&=\Inv(\s\1)\setminus\{(k,n)\mid k\in[n-1]\};\\
        \Neg(\r_\s\1)&=\Neg(\s\1)\setminus\{n\};\\
        \Nsp(\r_\s\1)&=\Nsp(\s\1)\setminus\{\{k,n\}\mid k\in[n-1]\}.
    \end{align*}
\end{lemma}

\begin{example}
    Consider $\s=3\o{6}15\o{4}2\in S^B_6$. Then $\s\1=361\o{5}4\o{2}$, so we can compute 
    \begin{eqnarray*}
    \Inv(\s\1) &=&\{(1,3),(1,4),(1,6),(2,3),(2,4),(2,5),(2,6),(3,4),(3,6),(5,6)\},\\
    \Neg(\s\1) &=&\{4,6\},\\
    \Nsp(\s\1) &=&\{\{1,4\},\{3,4\},\{3,6\},\{4,5\},\{4,6\}\}.
    \end{eqnarray*}
    Moreover, $\r_\s=315\o{4}2$ and $\r_\s\1=251\o{4}3$. Therefore, as stated in Lemma~\ref{lem:ncutInvNegNsp}, we have 
    \begin{eqnarray*}
    \Inv(\r_\s\1)&=&\{(1,3),(1,4),(2,3),(2,4),(2,5),(3,4)\},\\ 
    \Neg(\r_\s\1)&=&\{4\},\\
    \Nsp(\r_\s\1)&=&\{\{1,4\},\{3,4\},\{4,5\}\}.
    \end{eqnarray*}
\end{example} 

For simplicity, we identify strings with their underlying sets of elements. In particular, given two strings $B=b_1b_2\c b_k$ and $B'=b'_1b'_2\c b'_h$, we write $B\subseteq B'$ meaning that $\{b_1,\ldots,b_k\}\subseteq\{b'_1,\ldots ,b'_h\}$. Similarly, $B\cup B'$ denotes the set $\{b_1,\ldots,b_k\}\cup\{b'_1,\ldots ,b'_h\}$.

The next proposition highlights some properties of the $n$-cut of a signed permutation and the map $\r$, in connection with the weak order. 
These results will be crucial in the proof of Theorem~\ref{thm:join}.
\begin{proposition}\label{prop:invarianza}
    Let $\s,\q\in S^B_n$ and suppose $\s\w\q$ in $(S^B_n,\w)$; then the following hold:
    \begin{enumerate}
        \item $B^*_\s\subseteq B^*_\q$;
        \item $\r_\s\w\r_\q$ in $(S^B_{n-1},\w)$; that is, the map $\r$ is a surjective morphism of posets.
    \end{enumerate}
\end{proposition}

The converse implication of the first statement in Proposition~\ref{prop:invarianza} is not true in general.
Indeed, if we consider $\s=34\o{1}2$ and $\q=1\o{4}23$, with $4$-cuts $\bncut{\s}=(\o{2}1\o{3}3,\o{1}2)$ and $\bncut{\q}=(\o{2}\o{3},\o{1}123)$, then  $B^*_\s\subseteq B^*_\q$, but $\s \not\leq_R \q$, since, for example, $1\in\Neg(\s\1)\setminus\Neg(\q\1)$.

In the following example we collect some computations relative to multiple applications of the map $\r$. It is clear that if $\s \in S^B_n$ then $\r^i_\s \in S^B_{n-i}$. The symbol $\varepsilon$ denotes the empty string.
\begin{example}
    Consider $\s=\o{4}\o{1}5\o{2}3,\q=\o{4}5\o{1}3\o{2}\in S^B_5$; one can check that $\s\w\q$, for instance using \eqref{eq:weakstat}. It is not hard to verify the validity of the statements of Proposition~\ref{prop:invarianza} in each row of the table.
    \begin{center}
        \begin{tabular}{|c|c|c|c|c|}
            \hline
            & & & & \\[-4mm]
            $i$ & $\r^{i}_\s$ & $\r^i_\q$ & $\bncut{\r^{i}_\s}$ & $\bncut{\r^{i}_\q}$ \\[1.5mm]
            \hline
            & & & & \\[-4mm]
            $0$ & $\o{3}2\o{5}14\o{4}\o{1}5\o{2}3$ & $2\o{3}1\o{5}4\o{4}5\o{1}3\o{2}$ & $(\o{3}214\o{4}\o{1},\o{2}3)$ & $(2\o{3}14\o{4},\o{1}3\o{2})$ \\[0.5mm]
            \hline
            & & & & \\[-4mm]
            $1$ & $\o{3}214\o{4}\o{1}\o{2}3$ & $2\o{3}14\o{4}\o{1}3\o{2}$ & $(\o{3}21,\o{1}\o{2}3)$ & $(2\o{3}1,\o{1}3\o{2})$ \\[0.5mm]
            \hline
            & & & &\\[-4mm]
            $2$ & $\o{3}21\o{1}\o{2}3$ & $2\o{3}1\o{1}3\o{2}$ & $(21\o{1}\o{2},\varepsilon)$ & $(21\o{1},\o{2})$ \\[0.5mm]
            \hline
            & & & &\\[-4mm]
            $3$ & $21\o{1}\o{2}$ & $21\o{1}\o{2}$ & $(\varepsilon,1\o{1})$ & $(\varepsilon,1\o{1})$ \\[0.5mm]
            \hline
            & & & &\\[-4mm]
            $4$ & $1\o{1}$ & $1\o{1}$ & & \\[0.5mm]
            \hline
        \end{tabular}
    \end{center}
\end{example}

We now describe an algorithm that, given $\s,\q\in S^B_n$, produces another signed permutation $J(\s,\q)$. In the following description, we use the complete notation for signed permutations.

Given $\s\in S^B_n$, for any $2\leq i\leq n$, we define $$H^\s_i = \{a\in[\pm (i-1)]\mid \s\1(a)>\s\1(i)\}.$$
Observe that the underlying set of the string $B^*_{\r_\s^{n-i}}$ coincides with $H^\s_i$.

\begin{algorithm}[H]
    \caption{}\label{alg:algorithm}
    Consider $\s,\q\in S^B_n$.
    \begin{enumerate}
        \item Set 
        $$\pi_1=\begin{cases}
            \o{1}1,\quad&\emph{if }\s\1(1)>0\text{ and }\q\1(1)>0,\\
            1\o{1},&\emph{otherwise}.
        \end{cases}$$        
        \item For each $i\in\{2,3,\ldots,n\}$: 
        \begin{enumerate}
            \item Compute $H^\s_i$ and $H^\q_i$;
            \item Define $C^*_{i-1}$ as the shortest suffix of $\pi_{i-1}$ such that $$H^\s_i\cup H^\q_i\subseteq C^*_{i-1}$$
            and call $C_{i-1}$ the string such that $\pi_{i-1}=C_{i-1}C^*_{i-1}$.
            \item Define $\pi_i$ by placing $i$ in the complete notation of $\pi_{i-1}$ (consequently, $\o{i}$ is placed in the opposite entry) as follows:
            \begin{itemize}
                \item if $|C^*_{i-1}|=0$, place $i$ on the right of the last entry of $\pi_{i-1}$;
                \item if $|C^*_{i-1}|=i-1$, place $i$ immediately on the left of $C_{i-1}^*$ in case $\s\1(i)>0$ and $\q\1(i)>0$; otherwise, place $i$ immediately on the right of $C_{i-1}$;
                \item if $|C^*_{i-1}|\neq0$ and $|C^*_{i-1}|\neq i-1$, place $i$ immediately on the left of the first entry of $C^*_{i-1}$.
            \end{itemize}
        \end{enumerate}
        \item Set $J(\s,\q):=\pi_n$.
    \end{enumerate}
\end{algorithm}

\begin{remark}\label{rem:propCi}
    The definition of $\pi_i$ in Algorithm~\ref{alg:algorithm} can be visualized as follows:
        \begin{equation*}
            \pi_i=\begin{cases}
                C_{i-1}iC_{i-1}^{*'}\o{i}C^{*''}_{i-1},\ &\emph{if }|C_{i-1}^*|\geq i,\\
                C'_{i-1}\o{i}C''_{i-1}iC^*_{i-1},&\emph{if }|C_{i-1}^*|\leq i-2,\\
                C_{i-1}\o{i}iC_{i-1}^*,&\emph{if }|C_{i-1}^*|=i-1\emph{ and }\s\1(i),\q\1(i)>0,\\
                C_{i-1}i\o{i}C_{i-1}^*,&\emph{otherwise},
            \end{cases}
        \end{equation*}
        where 
        $$C_{i-1}=C'_{i-1}C''_{i-1}\emph{ with }|C_{i-1}'|=|C_{i-1}^*|$$
        and  
        $$C_{i-1}^*=C_{i-1}^{*'}C^{*''}_{i-1}\emph{ with }|C_{i-1}^{*''}|=|C_{i-1}|.$$
        Moreover, we observe that the $i$-cut of $\pi_i$ verifies $\bncut{\pi_i}=(C_{i-1},C_{i-1}^*)$, so $\r_{\pi_i}=\pi_{i-1}$.
\end{remark}

\begin{example}
    Let $\s=435\o{1}\o{2}, \q=\o{2}135\o{4}\in S^B_5$. We apply Algorithm~\ref{alg:algorithm}. Since $\s\1(1)<0$, we have $\pi_1=1\o{1}$. The following table displays the subsequent computations for $i$ from $2$ to $5$.
    \begin{center}
        \begin{tabular}{|c|c|c|c|c|c|}
            \hline
            & & & & &\\[-4mm]
            $i$ & $H^\s_i $ & $H^\q_i$ & $H^\s_i\cup H^\q_i$ & $C^*_{i-1}$ & $\pi_i$\\[1.5mm]
             \hline
            & & & & &\\[-4mm]
            $2$ & $\{\o{1},1\}$ &$\{1\}$ &$\{\o{1},1\}$ & $1\o{1}$ & $21\o{1}\o{2}$\\[0.5mm]
            \hline
            & & & & &\\[-4mm]
            $3$ &$\{\o{2},\o{1}\}$ & $\emptyset$ & $\{\o{2},\o{1}\}$ & $\o{1}\o{2}$  & $21\o{3}3\o{1}\o{2}$ \\[0.5mm]
            \hline
            & & & & &\\[-4mm]
            $4$ & $\{\o{1},\o{2},3\}$ & $\{\o{3},\o{2},\o{1},1,2,3\}$ & $\{\o{3},\o{2},\o{1},1,2,3\}$ & $21\o{3}3\o{1}\o{2}$  & $421\o{3}3\o{1}\o{2}\o{4}$ \\[0.5mm]
            \hline
            & & & & &\\[-4mm]
            $5$ & $\{\o{2},\o{1}\}$ & $\{\o{4}\}$ & $\{\o{4},\o{2},\o{1}\}$ & $\o{1}\o{2}\o{4}$ & $421\o{5}\o{3}35\o{1}\o{2}\o{4}$ \\[0.5mm]
            \hline
        \end{tabular}
    \end{center}
    To better understand the notation used in Remark~\ref{rem:propCi}, we note that in the previous computations we find:

    \begin{itemize}
        \item $C^{*'}_1=1\o{1}\quad \mbox{and}\quad C^{*''}_1=\varepsilon$;
        \item $C^{*'}_2=21\o{1}\o{2}\quad \mbox{and}\quad C^{*''}_2=\varepsilon$;
        \item $C^{*'}_3=21\o{3}3\o{1}\o{2}\quad \mbox{and}\quad C^{*''}_3=\varepsilon$;
        \item $C'_4=421\quad \mbox{and}\quad C''_4=\o{3}3$.
    \end{itemize}
    In conclusion, we obtain $J(\s,\q)=\pi_5=35\o{1}\o{2}\o{4}$.
\end{example}

The main result of this extended abstract is the following theorem, which shows that Algorithm~\ref{alg:algorithm} correctly computes the join of two signed permutations in $(S^B_n,\w)$.

\begin{theorem}\label{thm:join}
    Let $\s,\q\in S^B_n$; then 
    \begin{equation}
        J(\s,\q)=\s\jo\q.
    \end{equation}
\end{theorem}
We conclude this section by sketching our proof of Theorem~\ref{thm:join}. We prove by induction on $i$ that $\pi_i=\r_\s^{n-i}\jo\r_\q^{n-i}\in S^B_i$; in particular, we show that
    \begin{align*}
        \Inv((\r_\s^{n-i})\1)\cup\Inv((\r_\q^{n-i})\1)&\subseteq\Inv(\pi_i\1),\\
        \Neg((\r_\s^{n-i})\1)\cup\Neg((\r_\q^{n-i})\1)&\subseteq\Neg(\pi_i\1),\\
        \Nsp((\r_\s^{n-i})\1)\cup\Nsp((\r_\q^{n-i})\1)&\subseteq\Nsp(\pi_i\1),
    \end{align*}
getting  $\r_\s^{n-i},\r_\q^{n-i}\w\pi_i$, which implies $\r_\s^{n-i}\jo\r_\q^{n-i}\leq_R \pi_{i}.$
Then, in order to obtain $\pi_i\w\r_\s^{n-i}\jo\r_\q^{n-i}$, it is sufficient to prove that for any $\b\in S^B_i$ such that $\r_\s^{n-i},\r_\q^{n-i}\w\b$, one has $\pi_i\w\b$. This is done by using Lemma~\ref{lem:ncutInvNegNsp} and Proposition~\ref{prop:invarianza} on $\r_\s^{n-i},\r_\q^{n-i}\w\b$.
The case $i=n$, yields $\pi_n=J(\s,\q)=\s\jo\q.$

\section{Final remarks}

The main importance of Theorem~\ref{thm:join} and of its proof is that they provide an explicit description of the inversions, negative entries, and negative sum pairs of the join $\s\jo\q$ in terms of those of $\s$ and $\q$. 
This characterization reduces the proof of the following result to three cases and yields the inclusion $T_L(\s\jo\q)\subseteq \overline{T_L(\s)\cup T_L(\q)}$.

\begin{theorem}
    Let $\s,\q\in S^B_n$ and $t\in T_L(\s\jo\q)$. Then there exists a palindromic $(\s,\q)$-Bruhat path reaching $t$. In particular, $T_L(\s\jo\q)\subseteq \overline{T_L(\s)\cup T_L(\q)}$.
\end{theorem}

The reverse inclusion follows from \cite[Theorem 4.21]{dermenjian2025} of Dermenjian. Together, the two inclusions establish the validity of Conjecture~\ref{conjecture} for Coxeter groups of type $B$.

We conclude this extended abstract with a related open problem. As noted above, our algorithm generalizes the one defined by Markowsky for computing the join of two permutations in the weak order of the symmetric group $S_n$.
Indeed, $S_n$ is isomorphic to the subgroup of $S^B_n$ consisting of signed permutations with no negative entries.
Furthermore, $(S_n,\w)$ is a sublattice of $(S^B_n,\w)$ and Algorithm~\ref{alg:algorithm}, when restricted to this subgroup, coincides exactly with Markowsky's algorithm.

Within $S^B_n$, one can also consider $S^D_n$, the subgroup of index 2 consisting of all signed permutations with an even number of negative entries, see~\cite[Chapter 8]{bjorner2005combinatorics}. It is then natural to ask whether the same approach can be used to compute joins in type $D$.
In general, however, Algorithm~\ref{alg:algorithm} does not compute the join of two signed permutations in the poset $(S^D_n,\leq_D)$.

As a Coxeter system, $S^D_n$ can be generated by $S_D=\{s^D_0=\t{1}{\o{2}},s_1=\t{1}{2},\ldots,s_{n-1}=\t{n-1}{n}\}$ and the main difference with $S^B_n$ is the absence of reflections of the form $\t{a}{\o{a}}$. We denote by $\leq_D$ the (right) weak order defined on the Coxeter system $(S^D_n,S_D)$ and by $\leq_B$ the one induced by $S^B_n$ on the subgroup $S^D_n$. 
The two orders do not coincide. More precisely, given $\s,\q\in S^D_n$ one has $\s\leq_B \q$ implies $\s\leq_D \q$, but the converse implication is not true. One reason for this is that the left-reflection set of $\s\in S^D_n$ is given by $\left\{ \t{i}{j} \mid  (i,j) \in \Inv(\s\1)  \right\}\sqcup \left\{ \t{i}{\o{j}} \mid  \{i,j\} \in \Nsp(\s\1) \right\}$, which is a subset of its left-reflection set in $S^B_n$ described in \eqref{eq:leftref}.
The following example shows that Algorithm~\ref{alg:algorithm}, when applied to elements of $S^D_n$, may fail to produce their joins in $(S^D_n,\leq_D)$.

\begin{example}
    Consider the two signed permutations $\s=2314$ and $\q=\o{1}3\o{2}4$ in $S^D_4$. Their join in $(S^D_4,\leq_D)$ is $\s\vee_D\q=\o{3}1\o{2}4$, but if we apply Algorithm~\ref{alg:algorithm} we get $\pi_1=1\o{1}$ and 
  
    \begin{center}
        \begin{tabular}{|c|c|c|c|}
            \hline
            & & &\\[-4mm]
            $i$ & $ H^\s_i$& $H^\q_i$ & $\pi_i$ \\[1.5mm]
            \hline
            & & &\\[-4mm]
            $2$ & $\{1\}$ & $\{\o{1},1\}$ & $21\o{1}\o{2}$\\[0.5mm]
            \hline
            & & &\\[-4mm]
            $3$ & $\{1\}$ & $\{\o{2}\}$ & $231\o{1}\o{3}\o{2}$\\[0.5mm]
            \hline
            & & &\\[-4mm]
            $4$ & $\emptyset $ & $\emptyset$ & $\o{4}231\o{1}\o{3}\o{2}4$\\[0.5mm]
            \hline
        \end{tabular}
    \end{center}
    So, the join of $\s$ and $\q$ in $S^B_4$ is $\s\vee_B\q=\o{1}\o{3}\o{2}4$ which is not an element of $S^D_4$. 
\end{example}

At present, an analog to Algorithm~\ref{alg:algorithm} that operates directly on the complete notation of $S^D_n$ remains elusive, and
we propose its construction as an open problem. Such an algorithm would likely be a key instrument for extending our approach from \cite{BiagioliPerroneA} and obtaining a complete proof of Dyer's conjecture for type $D_n$.

\bibliographystyle{eptcs}
\bibliography{references}
\end{document}